\documentclass[UTF8,12pt]{article}
\usepackage{amsmath}
\usepackage{amsthm}
\usepackage[colorlinks, linkcolor=blue, anchorcolor=blue, citecolor=blue]{hyperref}
\usepackage{graphicx}
\usepackage{textcomp}
\usepackage{listings}
\usepackage{amsfonts}
\usepackage[textwidth=14.5cm]{geometry}
\usepackage{blindtext}
\newtheorem{thm}{Theorem}[section]
\newtheorem{cor}[thm]{Corollary}
\newtheorem{lemma}[thm]{Lemma}
\theoremstyle{definition}

\newtheorem{re}[thm]{Remark}
\title{Gradient estimate for $p$-heat equation on complete Riemannian manifold}
\author{Weihao Zheng}
\date{ }
\begin{document}
\maketitle
\begin{abstract}
We derive a Hamilton-type gradient estimate for positive solutions to the $p$-heat equation via the Moser iteration method. Our result relaxes the curvature assumption from sectional curvature to Ricci curvature. In addition, we obtain a Liouville-type theorem and a Harnack-type inequality as corollaries.
\end{abstract}
\section{Introduction}
The gradient estimate is an important tool in geometric analysis and partial differential equations. For example, it is widely used to derive Liouville-type theorems and Harnack inequalities, and to study eigenvalue problems. Cheng and Yau \cite{ChengYau75} obtained a local gradient estimate for positive harmonic functions on Riemannian manifolds. In 1986, Li and Yau studied this problem for the heat equation with a Schrödinger operator in their celebrated paper \cite{LiYau86}, where they derived the well-known Li-Yau-type gradient estimate. 

These methods have been generalized to various settings, for example, to the study of elliptic and parabolic equations \cite{MaZhu21,WuJY10,ZhangQ25}, and to weighted Riemannian manifolds \cite{LiXD05,WuJY10-2}. Moreover, other types of gradient estimates have been derived by choosing different auxiliary functions, such as the Hamilton-type \cite{Hamilton93}, the Souplet-Zhang-type \cite{SoupletZhang06}, the Li-Xu-type \cite{LiXu11}, and the Brighton-type \cite{Brighton13}.

In this paper we are interested in the $p$-Laplace operator. It is a natural generalization of the Laplace operator from the variational point of view. In 2009, Kotschwar and Ni \cite{KotschwarNi09} obtained local and global gradient estimates on Riemannian manifolds. More precisely, for  solutions to a parabolic equation, they established a local Li-Yau-type estimate. Throughout this paper, $B_o(R)$ denotes the geodesic ball of radius $R$ centered at a fixed point $o$, and we assume $B_o(R)\cap\partial M=\emptyset$. We state their result as follows. 
\begin{thm}\label{KN}
	{\rm(\cite{KotschwarNi09})} Let $(M^n,g)$ be a complete Riemannian manifold. Assume the sectional curvature $K_M$ satisfies $K_M\geq-K^2$ on $B_o(R)$. Suppose $u$ is a solution to
	\begin{equation*}
		u_t=\Delta_pu-\lvert\nabla u\rvert^p,
	\end{equation*}
	on $B_o(R)\times[0,T]$. If $1<p<2$ and $u_t\leq 0$ on $B_o(R)\times[0,T]$, then for any $\alpha>1$, there exist two constants $C_1=C_1(n,p,\alpha,KR)$ and $C_2=C_2(n,p)$ such that
	\begin{equation*}
		\lvert\nabla u\rvert^p+\alpha u_t\leq\max\bigl\{\big(8n\alpha^2(\frac{C_1}{R^2}+\frac{C_2}{t^{\frac{2}{p}}})\big)^{\frac{p}{2}},\ (2pn(n-1)K^2)^{\frac{p}{2}}\big(\frac{\alpha}{\alpha-1}\big)^p\bigr\},
	\end{equation*}
	for any $(x,t)\in B_o(\frac{R}{2})\times[0,T]$.
\end{thm}
 The $p$-Laplace operator is widely used in the study of eigenvalue problems \cite{Valtorta12,WangLi16}, the inverse mean curvature flow \cite{KotschwarNi09,Moser07}, $L^p$-logarithmic Sobolev inequalities \cite{WangWang23,WangYan21} and entropy formulas \cite{DelDolbealault03,WangYZ17}.

The method used by Kotschwar and Ni is mainly based on the maximum principle, which involves the second-order partial derivatives of the distance function. This necessitates the use of the Hessian comparison theorem and hence a condition on the sectional curvature, even for elliptic equations with the $p$-Laplace operator (see \cite{WangLF13,ZhaoShen20} for example).

In 2011, Wang and Zhang \cite{WangZhang11} introduced a new method based on the Moser iteration to obtain local gradient estimates for the $p$-harmonic equation, and established a Cheng-Yau-type estimate. Their method only differentiates the distance function once, hence no comparison theorem is needed, and the estimate holds with the Ricci curvature condition. Recently, the Moser iteration method has been widely used to derive gradient estimates for elliptic equations with the $p$-Laplace operator \cite{HeWangWei24,ZhaoYang18} and for parabolic equations with the Laplace operator \cite{HuangShen25,XiaQL23,ZhengWH26}.

In this paper, we derive a Hamilton-type gradient estimate on complete Riemannian manifolds under a Ricci curvature condition.
\begin{thm}\label{th1.1}
	Let $(M^n,g)$ be a complete Riemannian manifold whose Ricci curvature ${\rm Ric}$ is bounded from below, that is, ${\rm Ric}\geq-K,\ K\geq0$. Assume $u$ is a positive solution to
	\begin{equation}\label{1.1}
		u^{p-2}u_t=\Delta_pu,
	\end{equation}
	on $Q_{R,T}=B_o(R)\times(t_0-T,t_0]$, where $p>1$, and $0<D_1\leq u\leq D_2<\infty$ for two constants $D_1,D_2$. Then there exists a constant $C=C(n,p,\frac{D_2}{D_1})$ such that
	\begin{equation*}
		\frac{\lvert\nabla u\rvert}{u}(x,t)\leq C\big(\frac{1}{R}+\frac{1}{T^{\frac{1}{p}}}+\sqrt{K}\big),
	\end{equation*}
	for all $(x,t)\in B_o(\frac{R}{2})\times (t_0-\frac{T}{2},t_0]$.
\end{thm}
Compared with the result in \cite{KotschwarNi09}, our estimate requires no restriction on the time derivative, and the curvature condition is relaxed.
\begin{re}
	We also consider estimates involving both time and spatial derivatives, such as Li-Yau-type estimates. The main obstacle seems to be the non-homogeneity between the time derivative operator and the (linearized) $p$-Laplace operator.
\end{re}
As a corollary, we obtain the following Liouville-type theorem for ancient solutions to the $p$-heat equation.
\begin{cor}\label{cor1.2}
	Let $(M^n,g)$ be a complete noncompact Riemannian manifold without boundary whose Ricci curvature is nonnegative. Assume $u$ is a positive solution to (\ref{1.1}) on $M\times(-\infty,t_0]$ such that
	\begin{equation*}
		0<D_1\leq u\leq D_2<\infty,
	\end{equation*} 
	where $D_1,D_2$ are two constants. Then $u$ is constant.
\end{cor}  
Another corollary is the following Harnack-type inequality.
\begin{cor}\label{cor1.3}
	Let $(M^n,g)$ be a complete Riemannian manifold whose Ricci curvature ${\rm Ric}$ is bounded from below, that is, ${\rm Ric}\geq-K,\ K\geq0$. Assume $u$ is a positive solution to (\ref{1.1}) on $Q_{R,T}$, satisfying $0<D_1\leq u\leq D_2<\infty$. Then for any $x_1,x_2\in B_o(\frac{R}{4})$ and $t\in (t_0-\frac{T}{2},t_0]$, we have
	\begin{equation*}
		u(x_2,t)\exp\bigl(-Cl(\frac{1}{R}+\frac{1}{T^{\frac{1}{p}}}+\sqrt{K})\bigr)\leq u(x_1,t)\leq u(x_2,t)\exp\bigl(Cl(\frac{1}{R}+\frac{1}{T^{\frac{1}{p}}}+\sqrt{K})\bigr),
	\end{equation*}
	where $C=C(n,p,\frac{D_2}{D_1})$, and $l$ is the length of a minimal geodesic connecting $x_1$ and $x_2$.
\end{cor}
The paper is organized as follows. In Section 2, we briefly review some basic concepts of the $p$-heat equation and its linearized operator, then introduce a Sobolev embedding theorem. We derive the gradient estimate via the Moser iteration method in Section 3. In Section 4, we give the proofs of the corollaries.
\section{Preliminaries}
Let $(M^n,g)$ be a complete Riemannian manifold and $\Omega\subset M$ be a compact domain. For $p>1$, the $p$-energy of a function $\phi$ on $\Omega$ is defined by
\begin{equation*}
	E_p(\Omega,\phi)=\frac{1}{p}\int_{\Omega}\lvert\nabla\phi\rvert^p{\rm dvol}_g,
\end{equation*}
where ${\rm dvol}_g$ is the volume form induced by $g$. We omit ${\rm dvol}_g$ in the rest of the article for convenience.

A function $\phi$ is $p$-harmonic if it is a critical point of $E_p(\Omega,\cdot)$ for any compact $\Omega$. From the Euler-Lagrange equation we can infer that $\phi$ satisfies
\begin{equation*}
	\Delta_p\phi={\rm div}(\lvert\nabla\phi\rvert^{p-2}\nabla \phi)=0.
\end{equation*}
This naturally leads us to the concept of the $p$-heat equation. A function $u$ is a solution to the $p$-heat equation if $u$ satisfies
\begin{equation*}
	\Delta_pu={\rm div}(\lvert\nabla u\rvert^{p-2}\nabla u)=u^{p-2}u_t.
\end{equation*}

Let $u$ be a positive solution to (\ref{1.1}) and set $v=\log u$. Then $v$ satisfies
\begin{equation}\label{2.1}
	v_t=\Delta_pv+(p-1)\lvert\nabla v\rvert^p.
\end{equation}
The linearization of $\Delta_p$ at $v$, denoted by $L$, is
\begin{equation}\label{2.2}
	L(\phi)={\rm div}(\lvert\nabla v\rvert^{p-2}\nabla \phi)+(p-2){\rm div}(\lvert\nabla v\rvert^{p-4}\langle\nabla \phi,\nabla v\rangle\nabla v).
\end{equation}
In terms of a local orthonormal frame $\{e_1,\ldots,e_n\}$, the norm and trace with respect to $L$ are
\begin{equation*}
	\lVert \nabla^2 \phi\rVert^2_{L}=(g^{ij}+\frac{(p-2)v^iv^j}{\lvert\nabla v\rvert^2})(g^{kl}+\frac{(p-2)v^kv^l}{\lvert\nabla v\rvert^2})\phi_{ik}\phi_{jl},
\end{equation*}
and
\begin{equation*}
	Tr_L(\nabla^2\phi)=(g^{ij}+\frac{(p-2)v^iv^j}{\lvert\nabla v\rvert^2})\phi_{ij}.
\end{equation*}
In particular, for $v$ we have
\begin{equation}\label{2.3}
	\lVert \nabla^2 v\rVert^2_{L}=\lvert \nabla^2 v \rvert^2+\frac{p-2}{2\lvert\nabla v\rvert^2}\lvert\nabla\lvert\nabla v\rvert^2\rvert^2+\frac{(p-2)^2}{4\lvert\nabla v\rvert^4}\langle\nabla\lvert\nabla v\rvert^2,\nabla v\rangle^2,
\end{equation}
and
\begin{equation}\label{2.4}
	Tr_L(\nabla^2 v)=\Delta v+\frac{p-2}{2\lvert\nabla v\rvert^2}\langle \nabla\lvert\nabla v\rvert^2,\nabla v\rangle.
\end{equation}
\begin{re}
	Since the range of $p$ is $(1,\infty)$, the definitions above become degenerate on $\{\lvert \nabla v\rvert=0\}$. We follow the limiting approach adopted in \cite{HeWangWei24,WangZhang11} to handle this issue. Another approach is to use the $\epsilon$-regularization; see \cite{Hein07,WangLF18}.
\end{re}
At the end of this section we give a Sobolev embedding theorem which is essential for the Moser iteration method.
\begin{lemma}\label{le2.2}
	{\rm(\cite{SaloffCoste92})} Let $(M^n,g)$ be a complete Riemannian manifold with ${\rm Ric}\geq-K,\ K\geq0$, and let $B$ be a ball in $M$ with radius $R$ and volume $V$. For any $n>2$, there exists a constant $C(n)$ depending only on $n$ such that
	\begin{equation*}
		\Big(\int_B\lvert\phi\rvert^{\frac{2n}{n-2}}\Big)^{\frac{n-2}{n}}\leq e^{C(n)(1+\sqrt{K}R)}V^{-\frac{2}{n}}R^2\Big(\int_B(\lvert\nabla\phi\rvert^2+R^{-2}\phi^2)\Big),
	\end{equation*}
	for any $\phi\in C^{\infty}_0(B)$. For $n=2$, the above inequality still holds after replacing $n$ by any $n'>2$.
\end{lemma} 
\section{Proof of Theorem \ref{th1.1}}
We introduce the auxiliary functions first. Let
\begin{equation*}
	H=2+\frac{n(p-2)^2}{2(p-1)^2},\ h=H+\log D_2-v.
\end{equation*}
Thus, we get
\begin{equation*}
	H\leq h\leq H+\log(\frac{D_2}{D_1}),\ h_t=-v_t,\ \nabla h=-\nabla v.
\end{equation*}
We define
\begin{equation}\label{3.1}
	w=\frac{\lvert\nabla v\rvert^p}{h^p}.
\end{equation}

Direct computation gives the following lemma.
\begin{lemma}\label{le3.1}
	Let $v$ be a solution to (\ref{2.1}), and set $w=\frac{\lvert\nabla v\rvert^p}{h^p}$. Then we have
	\begin{equation*}\label{3.2}
		\begin{aligned}
			(L-\partial_t)w=&p(p-1)^2(h-1)h^{p-2}w^2\\
			&+p(p-1)(2-h)h^{p-3}w^{1-\frac{2}{p}}\langle\nabla w,\nabla v\rangle\\
			&+ph^{-p}\lvert\nabla v\rvert^{2p-4}\Big(\lVert \nabla^2v\rVert^2_L+(p-2)\frac{\lvert\nabla v\rvert^2}{h}Tr_L(\nabla^2v)\Big)\\
			&+ph^{-p}\lvert\nabla v\rvert^{2p-4}{\rm Ric}(\nabla v,\nabla v),
		\end{aligned}
	\end{equation*}
	on $\{\lvert\nabla v\rvert>0\}$.
\end{lemma}
\begin{proof}
	From the definition of $w$ we have $\lvert \nabla v\rvert^p=h^pw$. Taking the gradient on both sides gives
	\begin{equation}\label{3.3}
		\nabla w=\frac{p}{2}h^{-p}\lvert\nabla v\rvert^{p-2}\nabla\lvert \nabla v\rvert^2+ph^{-1}w\nabla v.
	\end{equation}
	Taking the scalar product with $\nabla v$, we obtain
	\begin{equation}\label{3.4}
		\langle\nabla \lvert\nabla v\rvert^2,\nabla v\rangle=\frac{2}{p}h^p\lvert\nabla v\rvert^{2-p}\langle\nabla w,\nabla v\rangle-2h^{p-1}\lvert\nabla v\rvert^{4-p}w.
	\end{equation}
	Applying the operator $L$ to the function $w$, we have
	\begin{equation*}
		Lw={\rm div}(\lvert\nabla v\rvert^{p-2}\nabla w)+(p-2){\rm div}(\lvert\nabla v\rvert^{p-4}\langle\nabla w,\nabla v\rangle\nabla v).
	\end{equation*}
	The first term on the right hand side is
	\begin{equation}\label{3.5}
		\begin{aligned}
			&{\rm div}(\lvert\nabla v\rvert^{p-2}\nabla w)\\
			=&{\rm div}(\frac{p}{2}h^{-p}\lvert\nabla v\rvert^{2p-4}\nabla\lvert\nabla v\rvert^2+ph^{-1}\lvert\nabla v\rvert^{p-2}w\nabla v)\\
			=&\frac{p^2}{2}h^{-p-1}\lvert\nabla v\rvert^{2p-4}\langle\nabla\lvert\nabla v\rvert^2,\nabla v\rangle+\frac{p}{2}h^{-p}\lvert\nabla v\rvert^{2p-4}\Delta\lvert\nabla v\rvert^2\\
			&+\frac{p(p-2)}{2}h^{-p}\lvert\nabla\lvert\nabla v\rvert^2\rvert^2+ph^{-2}\lvert\nabla v\rvert^pw\\
			&+\frac{p(p-2)}{2}h^{-1}\lvert\nabla v\rvert^{p-4}w\langle\nabla\lvert\nabla v\rvert^2,\nabla v\rangle\\
			&+ph^{-1}\lvert\nabla v\rvert^{p-2}\langle\nabla w,\nabla v\rangle+ph^{-1}\lvert\nabla v\rvert^{p-2}w\Delta v.
		\end{aligned}
	\end{equation}
	From (\ref{3.4}) we obtain
	\begin{equation}\label{3.6}
		\begin{aligned}
			&{\rm div}(\lvert\nabla v\rvert^{p-4}\langle\nabla w,\nabla v\rangle\nabla v)\\
			=&{\rm div}(\frac{p}{2}h^{-p}\lvert\nabla v\rvert^{2p-6}\langle\nabla\lvert\nabla v\rvert^2,\nabla v\rangle\nabla v+ph^{-1}\lvert\nabla v\rvert^{p-2}w\nabla v)\\
			=&\frac{p^2}{2}h^{-p-1}\lvert\nabla v\rvert^{2p-4}\langle\nabla\lvert\nabla v\rvert^2,\nabla v\rangle+\frac{p(p-3)}{2}h^{-p}\lvert\nabla v\rvert^{2p-8}\langle\nabla\lvert\nabla v\rvert^2,\nabla v\rangle^2\\
			&+\frac{p}{2}h^{-p}\lvert\nabla v\rvert^{2p-6}\langle\nabla\langle\nabla\lvert\nabla v\rvert^2,\nabla v\rangle,\nabla v\rangle+\frac{p}{2}h^{-p}\lvert\nabla v\rvert^{2p-6}\langle\nabla\lvert\nabla v\rvert^2,\nabla v\rangle\Delta v\\
			&+ph^{-2}\lvert\nabla v\rvert^pw+\frac{p(p-2)}{2}h^{-1}\lvert\nabla v\rvert^{p-4}w\langle \nabla \lvert\nabla v\rvert^2,\nabla v\rangle\\
			&+ph^{-1}\lvert\nabla v\rvert^{p-2}\langle\nabla w,\nabla v\rangle+ph^{-1}\lvert\nabla v\rvert^{p-2}w\Delta v.
		\end{aligned}
	\end{equation}
	By (\ref{2.1}), the time derivative satisfies 
	\begin{equation}\label{3.7}
		\begin{aligned}
			h^pw_t=&p\lvert\nabla v\rvert^{p-2}\langle\nabla v,\nabla v_t\rangle+ph^{p-1}wv_t\\
			=&p\lvert\nabla v\rvert^{2p-4}\langle \nabla v,\nabla \Delta v\rangle+\frac{p(p-2)}{2}\lvert\nabla v\rvert^{2p-6}\langle\nabla\lvert\nabla v\rvert^2,\nabla v\rangle\Delta v\\
			&+\frac{p(p-2)(p-4)}{4}\lvert\nabla v\rvert^{2p-8}\langle\nabla\lvert\nabla v\rvert^2,\nabla v\rangle^2\\
			&+\frac{p(p-2)}{2}\lvert\nabla v\rvert^{2p-6}\langle\nabla\langle\nabla\lvert\nabla v\rvert^2,\nabla v\rangle,\nabla v\rangle\\
			&+\frac{p^2(p-1)}{2}\lvert\nabla v\rvert^{2p-4}\langle\nabla\lvert\nabla v\rvert^2,\nabla v\rangle+ph^{p-1}\lvert\nabla v\rvert^{p-2}w\Delta v\\
			&+\frac{p(p-2)}{2}h^{p-1}\lvert\nabla v\rvert^{p-4}w\langle\nabla\lvert\nabla v\rvert^2,\nabla v\rangle+p(p-1)h^{p-1}\lvert\nabla v\rvert^pw.			
		\end{aligned}
	\end{equation}
	Combining (\ref{2.3}), (\ref{2.4}), (\ref{3.4})–(\ref{3.7}), and using the Bochner formula, we have
	\begin{equation}\label{3.8}
		\begin{aligned}
			(L-\partial_t)w=&p(p-1)^2(h-1)h^{p-2}w^2\\
			&+p(p-1)(2-h)h^{p-3}w^{1-\frac{2}{p}}\langle\nabla w,\nabla v\rangle\\
			&+ph^{-p}\lvert\nabla v\rvert^{2p-4}\Big(\lVert \nabla^2v\rVert^2_L+(p-2)\frac{\lvert\nabla v\rvert^2}{h}Tr_L(\nabla^2v)\Big)\\
			&+ph^{-p}\lvert\nabla v\rvert^{2p-4}{\rm Ric}(\nabla v,\nabla v).
		\end{aligned}
	\end{equation} 
\end{proof}
	We now continue the proof of Theorem \ref{th1.1}. The definition of $L$ yields
	\begin{equation*}
		\begin{aligned}
			&\lVert\nabla^2 v\rVert^2_L+(p-2)\frac{\lvert\nabla v\rvert^2}{h}Tr_L(\nabla^2v)\\
			\geq&\frac{1}{n}[Tr_L(\nabla^2v)+\frac{n(p-2)\lvert\nabla v\rvert^2}{2h}]^2-\frac{n(p-2)^2\lvert\nabla v\rvert^4}{4h^2}\\
			\geq&-\frac{n(p-2)^2\lvert\nabla v\rvert^4}{4h^2}.
		\end{aligned}
	\end{equation*} 
	Substituting this into (\ref{3.8}) and using the assumption of Ricci curvature and the definition of $w$, we obtain
	\begin{equation*}
		\begin{aligned}
			(L-\partial_t)w\geq& ph^{p-2}[(p-1)^2(h-1)-\frac{n(p-2)^2}{4}]w^2\\
			&+p(p-1)(2-h)h^{p-3}w^{1-\frac{2}{p}}\langle\nabla w,\nabla v\rangle\\
			&-pKh^{p-2}w^{2-\frac{2}{p}}.
		\end{aligned}
	\end{equation*}
	Since
	\begin{equation*}
		h\geq H=2+\frac{n(p-2)^2}{2(p-1)^2},
	\end{equation*}
	which implies
	\begin{equation*}
		(p-1)^2(h-1)-\frac{n(p-2)^2}{4}\geq\frac{1}{2}(p-1)^2h,
	\end{equation*}
	we have
	\begin{equation}\label{3.9}
		\begin{aligned}
			&w_t-Lw+\frac{p(p-1)^2}{2}h^{p-1}w^2\\
			\leq&-p(p-1)(2-h)h^{p-3}w^{1-\frac{2}{p}}\langle\nabla w,\nabla v\rangle\\
			&+pKh^{p-2}w^{2-\frac{2}{p}}.
		\end{aligned}
	\end{equation}
	We set $\alpha=\frac{p-2}{p}$, $w_\epsilon=\max\{w-\epsilon,0\}$, and introduce a cutoff function $\eta(x,t)$, which vanishes on the parabolic boundary of $Q_{R,T}$.
	
	We denote $Q_{r,t}=B_o(r)\times(t_0-t,t_0]$, $Q_{r,t,\tau}=B_o(r)\times(t_0-t,\tau]$. Multiplying both sides of (\ref{3.9}) by $\eta^2w_\epsilon^{q-1}$, where $q\geq2$ is a constant to be determined later, and integrating over $Q_{R,T,\tau}=B_o(R)\times(t_0-T,\tau]$, we obtain
	\begin{equation}\label{3.10}
		\begin{aligned}
			&\iint_{Q_{R,T,\tau}}\eta^2w_\epsilon^{q-1}w_t-\eta^2w_\epsilon^{q-1}Lw+\frac{p(p-1)^2}{2}h^{p-1}\eta^2w_\epsilon^{q-1}w^2\\
			\leq&\iint_{Q_{R,T,\tau}}-p(p-1)(2-h)h^{p-3}\eta^2w_\epsilon^{q-1}w^{1-\frac{2}{p}}\langle\nabla w,\nabla v\rangle\\
			&+\iint_{Q_{R,T,\tau}}pKh^{p-2}\eta^2w_\epsilon^{q-1}w^{2-\frac{2}{p}}.
		\end{aligned}
	\end{equation}
	The integration by parts formula yields
	\begin{equation}\label{3.11}
		\iint_{Q_{R,T,\tau}} \eta^2w_\epsilon^{q-1}w_t=\frac{1}{q}\int_{B_o(R)}\eta^2w^q_\epsilon\bigg|_{t=\tau}-\frac{2}{q}\iint_{Q_{R,T,\tau}}\eta\eta_tw^q_\epsilon.
	\end{equation}
	The second term on the left hand side of (\ref{3.10}) is
	\begin{equation}\label{3.12}
		\begin{aligned}
			-\iint_{Q_{R,T,\tau}}\eta^2w_\epsilon^{q-1}Lw=&(q-1)\iint_{Q_{R,T,\tau}}\Big(\eta^2w_\epsilon^{q-2}\lvert\nabla v\rvert^{p-2}\lvert\nabla w\rvert^2\\
			&+(p-2)\eta^2w_\epsilon^{q-2}\lvert\nabla v\rvert^{p-4}\langle\nabla w,\nabla v\rangle^2\Big)\\
			&+2\iint_{Q_{R,T,\tau}}\Big(\eta w_\epsilon^{q-1}\lvert\nabla v\rvert^{p-2}\langle\nabla w,\nabla \eta\rangle\\
			&+(p-2)\eta w_\epsilon^{q-1}\lvert\nabla v\rvert^{p-4}\langle\nabla w,\nabla v\rangle\langle\nabla \eta,\nabla v\rangle\Big).
		\end{aligned}
	\end{equation}
	Since 
	\begin{equation*}
		p>1,\ \lvert\nabla v\rvert^{p-2}=h^{p-2}w^{\alpha},\ H\leq h\leq H+\log(\frac{D_2}{D_1}),
	\end{equation*}
	there exists a constant $c=c(n,p,\frac{D_2}{D_1})$ such that
	\begin{equation}\label{3.13}
		\begin{aligned}
			c(q-1)\iint_{Q_{R,T,\tau}}\eta^2w^\alpha w_\epsilon^{q-2}\lvert\nabla w\rvert^2\leq&(q-1)\iint_{Q_{R,T,\tau}}\Big(\eta^2w_\epsilon^{q-2}\lvert\nabla v\rvert^{p-2}\lvert\nabla w\rvert^2\\
			&+(p-2)\eta^2w_\epsilon^{q-2}\lvert\nabla v\rvert^{p-4}\langle\nabla w,\nabla v\rangle^2\Big).
		\end{aligned}
	\end{equation}
	From now on we use $C$ or $c$ to denote a constant which depends only on $n,p,\frac{D_2}{D_1}$, and may differ from line to line. Similarly, we have
	\begin{equation}\label{3.14}
		\begin{aligned}
			&2\iint_{Q_{R,T,\tau}}\Big(\eta w_\epsilon^{q-1}\lvert\nabla v\rvert^{p-2}\langle\nabla w,\nabla \eta\rangle+(p-2)\eta w_\epsilon^{q-1}\lvert\nabla v\rvert^{p-4}\langle\nabla w,\nabla v\rangle\langle\nabla \eta,\nabla v\rangle\Big)\\
			\leq&C\iint_{Q_{R,T,\tau}}\eta w^\alpha w_\epsilon^{q-1}\lvert\nabla w\rvert\lvert\nabla \eta\rvert\\
			\leq&\frac{c(q-1)}{4}\iint_{Q_{R,T,\tau}}\eta^2w^\alpha w_\epsilon^{q-2}\lvert\nabla w\rvert^2+\frac{C}{q-1}\iint_{Q_{R,T,\tau}}w^\alpha w_\epsilon^q\lvert\nabla \eta\rvert^2.
		\end{aligned}
	\end{equation}
	Likewise,
	\begin{equation}\label{3.15}
		\begin{aligned}
			&\lvert\iint_{Q_{R,T,\tau}}-p(p-1)(2-h)h^{p-3}\eta^2w_\epsilon^{q-1}w^{1-\frac{2}{p}}\langle\nabla w,\nabla v\rangle\rvert\\
			\leq&C\iint_{Q_{R,T,\tau}}\eta^2w^{1-\frac{1}{p}}w_\epsilon^{q-1}\lvert\nabla w\rvert\\
			\leq& \frac{c(q-1)}{4}\iint_{Q_{R,T,\tau}}\eta^2w^\alpha w_\epsilon^{q-2}\lvert\nabla w\rvert^2+\frac{C}{q-1}\iint_{Q_{R,T,\tau}}\eta^2ww_\epsilon^q.
		\end{aligned}
	\end{equation}
	By the definition of $w_\epsilon$ we have $ww_\epsilon^q\leq w^2w_\epsilon^{q-1}$. Thus there exists a constant $q_{\min}$ which depends on $n,p,\frac{D_2}{D_1}$, such that for any $q\geq q_{\min}$, the last term of (\ref{3.15}) is absorbed. Combining (\ref{3.11})–(\ref{3.15}), after taking the supremum, we have
	\begin{equation}\label{3.16}
		\begin{aligned}
			&\frac{1}{q}\sup_{t_0-T\leq \tau\leq t_0}\int_{B_o(R)}\eta^2w^q_\epsilon+cq\iint_{Q_{R,T}}\eta^2w^\alpha w^{q-2}_\epsilon\lvert\nabla w\rvert^2+c\iint_{Q_{R,T}}\eta^2w^2w_\epsilon^{q-1}\\
			&\leq \frac{C}{q}\iint_{Q_{R,T}}\lvert\eta\eta_t\rvert w^q_\epsilon+\frac{C}{q}\iint_{Q_{R,T}}\lvert\nabla\eta\rvert^2w^\alpha w^q_\epsilon+CK\iint_{Q_{R,T}}\eta^2w^{1+\alpha}w^{q-1}_\epsilon,
		\end{aligned}
	\end{equation}
	since all the integrands on the right hand side are nonnegative.
	Letting $\epsilon\to 0$ and applying Fatou’s lemma, we deduce from (\ref{3.16}) that
	\begin{equation}\label{3.17}
		\begin{aligned}
			&\frac{1}{q}\sup_{t_0-T\leq \tau\leq t_0}\int_{B_o(R)}\eta^2w^q+cq\iint_{Q_{R,T}}\eta^2 w^{q-2+\alpha}\lvert\nabla w\rvert^2+c\iint_{Q_{R,T}}\eta^2w^{q+1}\\
			&\leq \frac{C}{q}\iint_{Q_{R,T}}\lvert\eta\eta_t\rvert w^q+\frac{C}{q}\iint_{Q_{R,T}}\lvert\nabla\eta\rvert^2w^{q+\alpha}+CK\iint_{Q_{R,T}}\eta^2w^{q+\alpha}.
		\end{aligned}
	\end{equation}
	For any $(x',t')\in B_o(\frac{R}{2})\times(t_0-\frac{T}{2},t_0]$, there exists a constant $c_0$ such that $\rho=c_0\min\{R,T^{\frac{1}{p}},K^{-\frac{1}{2}}\}$ satisfies
	\begin{equation*}
		B_{x'}(\rho)\times(t'-\rho^p,t']\subset B_o(R)\times(t_0-T,t_0],
	\end{equation*}
	and
	\begin{equation*}
		K\rho^2\leq1.
	\end{equation*}
	We denote
	\begin{equation*}
		Q_a=B_{x'}(a\rho)\times(t'-(a\rho)^p,t'],\ \frac{1}{2}\leq a\leq 1,
	\end{equation*}
	and set
	\begin{equation*}
		V_\rho=\lvert B_{x'}(\rho)\rvert.
	\end{equation*}
	Now we use the Moser iteration method to prove
	\begin{equation*}
		\sup_{Q_{1/2}} w\leq C\rho^{-p}.
	\end{equation*}
	We establish the initial bound first.
	\begin{lemma}\label{le3.2}
		Let $(x',t')\in B_o(\frac{R}{2})\times(t_0-\frac{T}{2},t_0]$ and $\rho=c_0\min\{R,T^{\frac{1}{p}},K^{-\frac{1}{2}}\}$. For a fixed $q'\geq q_{\min}$ and $\beta_0=q'+1$, there exists a constant $C=C(n,p,\frac{D_2}{D_1})$ such that
		\begin{equation}\label{3.18}
			\Big(\frac{1}{V_\rho\rho^p}\iint_{Q_{3/4}}w^{\beta_0}\Big)^{\frac{1}{\beta_0}}\leq C\rho^{-p}.
		\end{equation}
	\end{lemma}
	\begin{proof}
		 Let $\psi$ be a cutoff function which satisfies 
		\begin{equation*}
			\begin{aligned}
				&0\leq\psi\leq1,\ {\rm supp}\psi\subset Q_1,\ \psi=1\ \text{on}\ Q_{\frac{3}{4}},\\
				&\lvert\nabla\psi\rvert\leq C\rho^{-1},\ \lvert\psi_t\rvert\leq C\rho^{-p}.
			\end{aligned}
		\end{equation*}
		Set $\eta^2=\psi^m$, where $m$ is a constant satisfying
		\begin{equation*}
			m\geq q'+1,\ m(1-\alpha)\geq2(q'+1). 
		\end{equation*}
		Since $K\rho^2\leq1$, from (\ref{3.17}) we have
		\begin{equation}\label{3.19}
			\begin{aligned}
				c\iint_{Q_{1}}\psi^mw^{q'+1}\leq&C\rho^{-p}\iint_{Q_{1}}\psi^{m-1}w^{q'}\\
				&+C\rho^{-2}\iint_{Q_{1}}\psi^{m-2}w^{q'+\alpha}.
			\end{aligned}
		\end{equation}
		Young's inequality yields
		\begin{equation*}
			C\rho^{-p}\psi^{m-1}w^{q'}\leq\frac{c}{4}\psi^mw^{q'+1}+C\rho^{-p(q'+1)},
		\end{equation*}
		and
		\begin{equation*}
			C\rho^{-2}\psi^{m-2}w^{q'+\alpha}\leq\frac{c}{4}\psi^mw^{q'+1}+C\rho^{-p(q'+1)}.
		\end{equation*}
		Combining them with (\ref{3.19}) and dividing both sides by $\frac{c}{2}$, we obtain
		\begin{equation*}
			\iint_{Q_{3/4}}w^{q'+1}\leq CV_\rho\rho^p\rho^{-p(q'+1)}.
		\end{equation*}
		Hence (\ref{3.18}) holds.
	\end{proof}
	Now we give the energy estimate for the iteration. Without loss of generality, we may assume $n>2$ (when $n=2$ the proof is similar after replacing $n$ by $n'=4$).  We define
	\begin{equation*}
		\theta=\frac{n+2}{n},\ \sigma=\max\{\alpha,0\},\  \delta=\alpha-\theta\sigma.
	\end{equation*}
	and then set
	\begin{equation*}
		a_j=\frac{1}{2}+\frac{2^{-j}}{4},\ j=0,1,2,\ldots,
	\end{equation*}
	\begin{equation*}
		\beta_{j+1}=\theta\beta_j+\delta,\ q_j=\beta_j-\sigma.
	\end{equation*} 
	Then we have $q_j\geq q_0=\beta_0-\sigma>q'\geq q_{\min}$, and $\beta_j\to+\infty$ provided $q'\geq q_{\min}$ is chosen large enough so that
	\begin{equation*}
		\beta_0>-\frac{\delta}{\theta-1}.
	\end{equation*}
	Let $\{\eta_j\}_{j=0}^\infty$ be a family of cutoff functions satisfying
	\begin{equation*}
		\eta_j=1\ \text{in}\ Q_{a_{j+1}},\ {\rm supp}\eta_j\subset Q_{a_j},
	\end{equation*}
	\begin{equation*}
		\lvert\eta_{j,t}\rvert\leq C4^j\rho^{-p},\ \lvert\nabla\eta_j\rvert^2\leq C4^j\rho^{-2}.
	\end{equation*}
	The basic inequality gives
	\begin{equation*}
		\lvert\nabla\big(\eta w^{\frac{q+\alpha}{2}}\big)\rvert^2\leq2\lvert\nabla\eta\rvert^2 w^{q+\alpha}+\frac{(q+\alpha)^2}{2}\eta^2w^{q+\alpha-2}\lvert\nabla w\rvert^2.
	\end{equation*}
	Setting $q=q_j$ and $\eta=\eta_j$, we can infer from (\ref{3.17}) that
	\begin{equation}\label{3.20}
		\begin{aligned}
			&\sup_\tau\int_{B_{x'}(a_j\rho)}\eta_j^2w^{q_j}+\iint_{Q_{a_j}}\lvert\nabla\big(\eta_jw^{\frac{(q_j+\alpha)}{2}}\big)\rvert^2+\rho^{-2}\iint_{Q_{a_j}}\eta_j^2w^{q_j+\alpha}\\
			\leq&Cq_j^24^j\Big(\rho^{-p}\iint_{Q_{a_j}}w^{q_j}+\rho^{-2}\iint_{Q_{a_j}}w^{q_j+\alpha}\Big).
		\end{aligned}
	\end{equation}
	Considering the cases $w\geq\rho^{-p}$ and $w<\rho^{-p}$ separately, we obtain
	\begin{equation*}
		w^{q_j}\leq \rho^{p\sigma}w^{\beta_j}+\rho^{-pq_j},
	\end{equation*}
	and
	\begin{equation*}
		w^{q_j+\alpha}\leq\rho^{p(\sigma-\alpha)}w^{\beta_j}+\rho^{-p(q_j+\alpha)}.
	\end{equation*}
	Using the fact that $p\alpha=p-2$, (\ref{3.20}) becomes
	\begin{equation}\label{3.21}
		\begin{aligned}
			&\sup_\tau\int_{B_{x'}(a_j\rho)}\eta_j^2w^{q_j}+\iint_{Q_{a_j}}\lvert\nabla\big(\eta_jw^{\frac{(q_j+\alpha)}{2}}\big)\rvert^2+\rho^{-2}\iint_{Q_{a_j}}\eta_j^2w^{q_j+\alpha}\\
			\leq&Cq_j^24^j\rho^{-p(1-\sigma)}\Big(\iint_{Q_{a_j}}w^{\beta_j}+V_{\rho}\rho^p\rho^{-p\beta_j}\Big).
		\end{aligned}
	\end{equation}
	For each fixed time level, Hölder's inequality gives
	\begin{equation*}
		\begin{aligned}
			&\int_{B_{x'}(a_j\rho)}\eta_j^{\frac{2n+4}{n}}w^{\frac{q_j(n+2)}{n}+\alpha}\\
			\leq&\Big(\int_{B_{x'}(a_j\rho)}\big(\eta_jw^{\frac{q_j+\alpha}{2}}\big)^{\frac{2n}{n-2}}\Big)^{\frac{n-2}{n}}\Big(\int_{B_{x'}(a_j\rho)}\eta_j^2w^{q_j}\Big)^\frac{2}{n}.
		\end{aligned}
	\end{equation*}
	Integrating over time and using Lemma \ref{le2.2}, we have
	\begin{equation}\label{3.22}
		\begin{aligned}
			\iint_{Q_{a_{j+1}}}w^{\theta q_j+\alpha}\leq&C\rho^2V_{\rho}^{-\frac{2}{n}}\Big(\sup_\tau\int_{B_{x'}(a_j\rho)}\eta_j^2w^{q_j}\Big)^{\frac{2}{n}}\\
			&\times\Big(\iint_{Q_{a_j}}\lvert\nabla\big(\eta_jw^{\frac{q_j+\alpha}{2}}\big)\rvert^2+\rho^{-2}\iint_{Q_{a_j}}\eta_j^2w^{q_j+\alpha}\Big).
		\end{aligned}
	\end{equation}
	Note that the volume $V_\rho$ is unchanged in each step of the iteration process, since all functions are supported in $B_{x'}(\rho)$. Combining (\ref{3.21}) and (\ref{3.22}) leads to
	\begin{equation}\label{3.23}
		\begin{aligned}
			\Big(\iint_{Q_{a_{j+1}}}w^{\theta q_j+\alpha}\Big)^{\frac{1}{\theta}}\leq&Cq_j^24^j\rho^{\frac{2}{\theta}-p(1-\sigma)}V_\rho^{-\frac{2}{n\theta}}\\
			&\times\Big(\iint_{Q_{a_j}}w^{\beta_j}+V_\rho
			\rho^p\rho^{-p\beta_j}\Big).
		\end{aligned}
	\end{equation}
	The definitions of $\delta,\sigma,\theta$ imply
	\begin{equation*}
		\theta q_j+\alpha=\theta\beta_j+\delta,
	\end{equation*}
	\begin{equation*}
		\frac{2}{n\theta}+\frac{1}{\theta}=1,
	\end{equation*}
	\begin{equation*}
		\frac{2}{\theta}-p(1-\sigma)+p(1-\frac{1}{\theta})=-\frac{p\delta}{\theta}.
	\end{equation*}
	For convenience we denote 
	\begin{equation*}
		N_{\beta_j}(Q_{a_j})=\Big(\frac{1}{V_\rho\rho^p}\iint_{Q_{a_j}}w^{\beta_j}+\rho^{-p\beta_j}\Big)^{\frac{1}{\beta_j}}.
	\end{equation*}
	Then from (\ref{3.23}) we obtain
	\begin{equation}\label{3.24}
		\Big(\frac{1}{V_\rho\rho^p}\iint_{Q_{a_{j+1}}}w^{\theta\beta_j+\delta}\Big)^{\frac{1}{\theta}}\leq Cq_j^24^j\rho^{-\frac{p\delta}{\theta}}N_{\beta_j}(Q_{a_j})^{\beta_j}.
	\end{equation}
	Moreover, $\rho^{-\frac{p(\theta\beta_j+\delta)}{\theta}}$ is also controlled by the right hand side of (\ref{3.24}) since
	\begin{equation*}
		\rho^{-p\beta_j}\leq N_{\beta_j}(Q_{a_j})^{\beta_j}.
	\end{equation*}
	Thus, from the fact that
	\begin{equation*}
		q_j=\beta_j-\sigma\leq\beta_j,
	\end{equation*}
	we have
	\begin{equation*}
		\begin{aligned}
			N_{\theta\beta_j+\delta}(Q_{a_{j+1}})\leq&(C\beta_j^24^j)^{\frac{\theta}{\theta\beta_j+\delta}}\rho^{-\frac{p\delta}{\theta\beta_j+\delta}}\\
			&\times N_{\beta_j}(Q_{a_j})^{\frac{\theta\beta_j}{\theta\beta_j+\delta}},
		\end{aligned}
	\end{equation*}
	that is,
	\begin{equation}\label{3.25}
		N_{\beta_{j+1}}(Q_{a_{j+1}})\leq(C\beta_j^24^j)^{\frac{\theta}{\beta_{j+1}}}\rho^{-\frac{p\delta}{\beta_{j+1}}}N_{\beta_j}(Q_{a_j})^{\frac{\theta\beta_j}{\beta_{j+1}}},
	\end{equation}
	for all $j=0,1,2,\ldots$.
	
	The iteration process yields
	\begin{equation}\label{3.26}
		\begin{aligned}
			N_{\beta_j}(Q_{a_j})\leq&\Big(\prod_{i=0}^{j-1}(C\beta_i^24^i)^{\frac{\theta^{j-i}}{\beta_j}}\Big)\rho^{\frac{p\theta^j\beta_0}{\beta_j}-p}N_{\beta_0}(Q_{\frac{3}{4}})^{\frac{\theta^j\beta_0}{\beta_j}}.
		\end{aligned}
	\end{equation}
	The products are uniformly bounded since
	\begin{equation*}
		\frac{\theta^{j-i}}{\beta_j}\leq\frac{\theta^{-i}}{\beta_0+\frac{\delta}{\theta-1}},\ \log(C\beta_i^24^i)\leq C(1+i),
	\end{equation*}
	and
	\begin{equation*}
		\sum_{i=0}^{\infty}(1+i)\theta^{-i}<\infty.
	\end{equation*}
	From Lemma \ref{le3.2} we have
	\begin{equation*}
		N_{\beta_0}(Q_{\frac{3}{4}})\leq C\rho^{-p}.
	\end{equation*}
	Then (\ref{3.26}) becomes
	\begin{equation*}
		\begin{aligned}
			N_{\beta_j}(Q_{a_j})\leq&C\rho^{\frac{p\theta^j\beta_0}{\beta_j}-p}(C\rho^{-p})^{\frac{\theta^j\beta_0}{\beta_j}}\\
			\leq&C\rho^{-p}.
		\end{aligned}
	\end{equation*}
	Letting $j\to\infty$, we obtain
	\begin{equation*}
		\sup_{Q_{1/2}} w\leq C\rho^{-p}.
	\end{equation*}
	Thus, for any $(x',t')\in B_o(\frac{R}{2})\times(t_0-\frac{T}{2},t_0]$, we have
	\begin{equation*}
		\frac{\lvert\nabla u\rvert}{u}=\lvert\nabla v\rvert=hw^{\frac{1}{p}}\leq C^{\frac{1}{p}}\rho^{-1}\cdot(H+\log(\frac{D_2}{D_1})).
	\end{equation*}
	The definition of $\rho$ gives 
	\begin{equation*}
		\rho^{-1}\leq C(R^{-1}+T^{-\frac{1}{p}}+\sqrt{K}).
	\end{equation*}
	Therefore,
	\begin{equation*}
		\sup_{B_o(R/2)\times(t_0-T/2,t_0]}\frac{\lvert\nabla u\rvert}{u}\leq C(n,p,\frac{D_2}{D_1})(R^{-1}+T^{-\frac{1}{p}}+\sqrt{K}).
	\end{equation*}
	This completes the proof of Theorem \ref{th1.1}.
\section{Proofs of the Corollaries}
\begin{proof}[Proof of Corollary \ref{cor1.2}]
	Since $M$ has nonnegative Ricci curvature and $u$ is a positive ancient solution, Theorem \ref{th1.1} yields
	\begin{equation*}
		\frac{\lvert\nabla u\rvert}{u}\leq C(n,p,\frac{D_2}{D_1})(\frac{1}{R}+\frac{1}{T^{\frac{1}{p}}}),
	\end{equation*}	
	for any $(x,t)\in B_o(\frac{R}{2})\times(t_0-\frac{T}{2},t_0]$ and $R,T>0$. Letting $R$ and $T$ tend to $\infty$, we obtain $\lvert\nabla u\rvert=0$. From (\ref{1.1}) we have $u_t=0$ as well, and the proof is complete.
\end{proof} 
\begin{proof}[Proof of Corollary \ref{cor1.3}]
	Let $\gamma(s)$ be a minimal geodesic connecting $x_1$ and $x_2$ with $\gamma(0)=x_1$ and $\gamma(l)=x_2$. Clearly the distance between $x_1$ and $x_2$ is less than $\frac{R}{2}$; hence $\gamma$ lies inside $B_o(\frac{R}{2})$. Using Theorem \ref{th1.1}, we obtain
	\begin{equation*}
		\begin{aligned}
			\log u(x_2,t)-\log u(x_1,t)=&\int_{0}^{l}\langle\nabla \log u,\gamma'(s)\rangle ds\\
			\leq&\int_{0}^l\frac{\lvert\nabla u\rvert}{u}ds\\
			\leq&\int_{0}^lC(n,p,\frac{D_2}{D_1})(R^{-1}+T^{-\frac{1}{p}}+\sqrt{K})ds\\
			\leq& Cl(R^{-1}+T^{-\frac{1}{p}}+\sqrt{K}).
		\end{aligned}
	\end{equation*}
	The reverse inequality follows similarly by interchanging $x_1$ and $x_2$.
\end{proof}
\section*{Declaration}
\subsection*{Competing interests}
The author declares no conflict of interests.
\subsection*{Data availability}
No data was generated or analyzed during the current study.
\subsection*{AI usage}
During the preparation of this work, the author used generative AI tools to assist with computation and verification, and to refine the language in the presentation of the manuscript. The technical methodology was conceived independently by the author. All calculations and logical derivations have been checked by the author, who takes full responsibility for the content of this article.
\bibliography{bib}
\bibliographystyle{plain}
Weihao Zheng

School of Mathematical Sciences, Xiamen University, 

Xiamen, 361005, PR China

\itshape{Email address}: whzheng@xmu.edu.cn
\end{document}